\documentstyle[11pt]{article}

\input amssym.def
\input amssym
\input xypic

\def\demo{\noindent{\bf Proof. }}
\def\QED{\hfill$\Box$}

\def\dim{{\rm dim }}
\def\ass{{\rm Ass }}
\newtheorem{Theorem}{Theorem}[section]

\newtheorem{Corollary}[Theorem]{Corollary}
\newtheorem{Proposition}[Theorem]{Proposition}
\newtheorem{Remark}[Theorem]{Remark}

\begin{document}
\title{
{{\Huge \bf Links of Prime Ideals}} \\
\footnotetext{AMS 1991 {\em Mathematics   Subject Classification}.
Primary 13H10; Secondary 13C40,  13D40, 13D45, 13H15.}
 }

\author{Alberto Corso \and
Claudia Polini\thanks{The author was partially supported by NSF grant
STC-91-19999.} \and
Wolmer V. Vasconcelos\thanks{The author was partially supported by
NSF grant DMS-93-43436.} \\
\begin{tabular}{c}
\vspace{0.2in} \      \\
Department of Mathematics, Rutgers University \\
New Brunswick, New Jersey 08903 \\
e-mail: corso, polini, vasconce\,@math.rutgers.edu
\end{tabular}}

\date{ }

\maketitle \vspace{-0.2in}

\begin{abstract}
\noindent We exhibit the elementary but somewhat surprising
property that most direct links
of prime ideals in Gorenstein rings are equimultiple ideals.
It leads to the construction of a bountiful set of Cohen--Macaulay
Rees algebras.
\end{abstract}

\section{Introduction}
{\it Grosso modo}, a link of an ideal ${\frak p}$ of a Noetherian ring
$R$ is an ideal of the form $I = ({\bf z})\colon {\frak p}$, where
${\bf z} = z_1, \ldots, z_g$ is a regular sequence and $g$ is the
codimension of ${\frak p}$. This is a very common operation in
commutative algebra, particularly in duality theory,
and plays an important role in current methods to
effect primary decomposition of polynomial ideals (see \cite{EHV}).

More precisely,
let $(R, {\frak m})$ be a Cohen--Macaulay local ring, and let $I$ and
$J$ be two ideals. $I$ and $J$ are said to be {\em linked} if there exists
a regular sequence ${\bf z}= z_1, \ldots, z_g\subset I\cap J$ such that
$I= ({\bf z})\colon J$ and $J= ({\bf z})\colon I$.
If $R$ is a Gorenstein ring and $J$ is an unmixed ideal, then for any
regular sequence ${ z_1, \ldots, z_g}\subset J$,
 where $g$ is the codimension of
$J$,  the
ideals $J$ and $I = ({\bf z})\colon J$ are linked (see \cite{pszpiro}). In
this case we simply say that $I$ is a {\em link} of $J$.

To frame our results we recall a notion of Northcott--Rees (see \cite{NR}).
An ideal $J\subset I$ is a {\em {reduction}} of $I$ if $JI^r = I^{r+1}$
for some integer $r$; the least such integer
 is the {\em {reduction number}}
of $I$ with respect to $J$. Phrased otherwise, this means that
$$
\diagram
R[Jt]= \sum\limits_{n=0}^{\infty } J^{n}t^{n} \rto|<<\tip &
\sum\limits_{n=0}^{\infty }
I^{n}t^{n}=R[It]
\enddiagram
$$
is a finite morphism of the associated Rees algebras, and $r$ is the
bound of the degrees required to generate $R[It]$ as a module over $R[Jt]$.
If $J$ is generated by a regular sequence
 the ideal $I$ is said to be {\em equimultiple}.

There are comparatively few general classes of ideals in
Cohen--Macaulay rings whose associated Rees algebra is
Cohen--Macaulay. The most basic case, of an ideal $I$ generated by a
regular sequence ${\bf z}= z_1, \ldots, z_g$, has a very explicit
presentation as a quotient of a polynomial ring
\[R[It] \simeq R[T_1, \ldots,
T_g]/I_2\left(\begin{array}{ccc} z_1 &\cdots & z_g \\ T_1 & \cdots & T_g
\end{array} \right),  \]
where the relations are the maximal minors of a matrix whose two rows
are essentially a set of indeterminates. From \cite{EN} it follows
that $R[It]$ is Cohen--Macaulay.

Other classes of ideals, with Cohen--Macaulay Rees algebras, include
maximal ideals with restrictions on their multiplicities, classical
determinantal ideals attached to generic matrices,  ideals whose
Koszul homology modules have conveniently high depths, and ideals
having nice reductions.


%

The ideals studied here are
 links of a  prime ideal ${\frak p}$.
Most of them turn out to be equimultiple,  of reduction number one,
provided the localization $R_{\frak p}$ is a Gorenstein ideal. As a
consequence, if $R$
is a Gorenstein ring and ${\frak p}$ is a Cohen--Macaulay ideal then
the corresponding associated graded ring will be Cohen--Macaulay.
Since ideals such as ${\frak p}$ are very common, the process described
here provides Cohen--Macaulay algebras with a great deal of profusion.

\section{Equimultiple Ideals}

The relationship between a regular sequence and the corresponding
 link of a prime is discussed next. The general
terminology used here is that of
\cite{Mat}; the basic reference for local duality and the theory of
the canonical module will be \cite{HK}.

 Two pieces of notation will be used: if $ R $ is a local ring
and $M$ is a finitely generated $ R$-module, $\nu(M)$ denotes its minimal
number of generators; if furthermore $M$ is Artinian, we denote its
length by $\lambda(M)$.

\begin{Theorem} \label{primelink1} Let $R$ be a Cohen--Macaulay ring,
  ${\frak p}$  a prime ideal of codimension $g$, and let ${\bf z}= z_1,
  \ldots, z_g\subset {\frak p}$ be a regular sequence. Set $J= ({\bf z})$
  and $I =J \colon {\frak p}$. Suppose that $R_{{\frak p}}$ is a Gorenstein ring.
  Then $I$ is an equimultiple ideal with reduction number one, more
  precisely,
$$
I^2=J I
$$
in the following two cases:
\begin{itemize}
\item[{\rm (a)}] $R_{{\frak p}}$ is not a regular local ring;
\item[{\rm (b)}] $R_{{\frak p}}$ is a regular local ring of dimension at
least $2$ and two elements in the sequence ${\bf z}$ lie in the
symbolic square ${\frak p}^{(2)}$.
\end{itemize}
\end{Theorem}

The proof requires the assemblage of 3 elementary techniques. We deal
with them separately first.

\subsection*{Reduction to maximal ideals}

To establish the equality $I^2= J I$, we examine the
associated prime ideals of $J I$. To that end, looking at the
exact sequence (see \cite{HH})
$$
\diagram
0 \rto & J/JI \rto & R/JI \rto &  R/J \rto & 0,
\enddiagram
$$
since \[J/JI = J/J^2 \otimes R/I = (R/I)^g,\]
it follows that the associated primes of $JI$ satisfy the condition
$$
\ass(R/JI) \subset \ass(R/I) \cup \ass(R/J).
$$
Therefore all the associated primes of $JI$ have codimension $g$, and
since $I$ and $J$ are equal at any localization $R_{\frak q}$, where
$\frak q$ is a prime of codimension $g$, distinct from ${\frak p}$, we may
finally assume that $R$ is a local ring and ${\frak p}$ its
maximal ideal.

\subsection*{Syzygetic ideals}

In making comparisons between the primary ideals $I^2$ and $JI$ we
consider the following diagram:

$$
\diagram
\ & I\drline^{\lambda(I/J)}\dlline_{\lambda(I/I^2)} &\  \\
I^2 \drline_{\lambda(I^2/JI)}&\  & J \dlline^{d\lambda(R/I)} \\
\ & JI &\
\enddiagram
$$

\vspace{1cm}

\noindent We need to establish some relationships between these lengths.
This will be done in the next result which contains our main
computation.

\begin{Proposition}\label{P1}
Let $(R, \frak m)$ be a Cohen--Macaulay local ring of
Krull dimension $d$,
let $J=(a_1, \ldots, a_d)$ be a $\frak m$--primary ideal
 and let $J \subset I$. Suppose $I$ is minimally
generated by  $d+r$ elements $(J, b_1, \ldots, b_r)$.
Denote
by $H_1$ the 1--dimensional Koszul homology module of $I$ and by
$\delta(I)$ the kernel of the natural surjection from the symmetric
square of $I$ onto $I^2$.
Then
\begin{equation}
\lambda(I/J) = \lambda(I^2/JI) + r \lambda(R/I) -\lambda(H_1)
+ \lambda(\delta(I)).\label{lengthS1}
\end{equation}
\end{Proposition}

\demo Let
$$
\diagram
0 \rto & Z \rto & {R}^{d+r}\rto & I \rto & 0,
\enddiagram
$$
be a minimal presentation of
$I$. Tensoring with $R/I$ we obtain the exact sequence (see \cite{SV3})
\begin{equation}
\diagram
0 \rto & \delta(I) \rto & H_1 \rto & (R/I)^{d+r}\rto & I/I^2 \rto & 0,
\enddiagram\label{syzygetic1}
\end{equation} where $\delta(I)$ accounts for the cycles in $Z$ whose
coefficients lie in $I$.

This
 leads to $
\lambda(I/I^2)-(d+r)\lambda(R/I)+\lambda(H_1)-\lambda(\delta(I))=0,
$ or equivalently
\begin{equation}
\lambda(I/I^2)-d\lambda(R/I)=r\lambda(R/I)-\lambda({H}_1)+\lambda(\delta(I)).
\label{eq1}
\end{equation}
To account for  the left-hand side of (\ref{eq1}),
 consider the following two exact sequences
$$
\diagram
0 \rto & J/JI \rto & I/JI\rto & I/J \rto & 0,
\enddiagram
$$
$$
\diagram
0 \rto & I^2/JI \rto & I/JI\rto & I/I^2 \rto & 0.
\enddiagram
$$
From the first one we get that
$ \lambda(I/JI)=\lambda(J/JI)+\lambda(I/J), $ while from the second one we get
that $ \lambda(I/JI)=\lambda(I^2/JI)+\lambda(I/I^2). $
By comparing the last two equations, we have that $
\lambda(J/JI)+\lambda(I/J)=\lambda(I^2/JI)+\lambda(I/I^2), $
or equivalently
\begin{equation}
\lambda(I/I^2)-\lambda(J/JI)=\lambda(I/J)-\lambda(I^2/JI).
\label{eq2}
\end{equation} As we remarked earlier,
 $ J/JI $ is a free $R/I$--module of rank $d$, so from
 $ \lambda(J/JI)=d\lambda(R/I) $,
the last equation  becomes:
\begin{equation}
\lambda(I/I^2)- d\lambda(R/I)=\lambda(I/J)-\lambda(I^2/JI).\label{eqx}
\end{equation}
If we compare (\ref{eqx}) with (\ref{eq1}), we get $
\lambda(I/J)-\lambda(I^2/JI)=r\lambda(R/I)-\lambda(H_1)+\lambda(\delta(I)), $ or equivalently
$$
\lambda(I/J)=\lambda(I^2/JI)+r\lambda(R/I)-\lambda(H_1)+\lambda(\delta(I)),
$$
as desired.\QED

\begin{Corollary}\label{C1}
Under the same assumption as in {\rm Proposition \ref{P1}},
if $R$ is a Gorenstein local ring and
 $r=1$ then
\begin{equation} \label{lengthS2}
\lambda(I/J) = \lambda(I^2/JI)+ \lambda(\delta(I)), \ {\rm \
with } \ \delta(I)\neq 0.
\end{equation}
\end{Corollary}

\demo
Since $ R $ is Gorenstein and $ R/I $ is Cohen-Macaulay,
the canonical module of $ R/I $ is the last non-vanishing
Koszul homology module
$$
\omega_{R/I}={\rm Ext}_R^d(R/I, R)={\rm Hom}_{R/J}(R/I,R/J)=(J\colon I)/J=H_1,
$$
as $ I $ is generated by a regular sequence
plus an extra element.
Recall that
 the canonical module of an Artinian ring is the injective envelope of
its residue field and therefore $ \lambda(H_1)=\lambda(R/I). $
If we apply Proposition \ref{P1}  we get
$$
\lambda(I/J)=\lambda(I^2/JI)+\lambda(\delta(I)).
$$
On the other hand, if $ \delta(I)=0, $ from (\ref{syzygetic1}),
we get the short exact sequence
$$
\diagram
0 \rto & H_1 \rto & (R/I)^{d+1}\rto & I/I^2 \rto & 0,
\enddiagram
$$
which splits since $ H_1 $ is injective. Therefore
$$
(R/I)^{d+1}\simeq H_1\oplus I/I^2.
$$
Now, $ \lambda(H_1)=\lambda(R/I) $ implies that $ H_1\simeq R/I $ and
so
$ I/I^2 \simeq  (R/I)^d, $  which implies that $ \nu(I)=d $ (since $ \nu(I)=\nu(I/I^2)
$ by Nakayama lemma), contradicting the assumption. \QED

\begin{Remark}{\rm
A similar computation would show the following more general formula.
Suppose that $(R,{\frak m})$ is a Cohen--Macaulay ring of type $s\geq
2$. Then
\[ \lambda(I^2/JI) + \lambda(\delta(I)) = {{s+1}\choose{2}}.\]
}\end{Remark}

\subsection*{Northcott ideals}

Finally, we recall an useful construction introduced by Northcott
(see \cite{Nor}), that provides for a very explicit description of a family
of links.

Let $R$ be a Noetherian ring,
${\bf u} = u_1, \ldots, u_n$  a sequence of elements in $R$,
and  $\varphi=(a_{ij})$  a $n \times n$ matrix with entries in
$R$. Denote by ${\bf v}= v_1, \ldots, v_n$  the sequence
\[{\bf [v]}^t= \varphi \cdot {\bf [u]}^t.\]
If  $\mbox{\rm grade } ({\bf v}) = n$, the ideal $N=({\bf v}, \det
\varphi )$ (where $\det \varphi$ is the determinant of the matrix
$\varphi$)  is known as the {\em Northcott ideal}
associated to the sequence\index{Northcott ideal} ${\bf u}$ and the
matrix $\varphi$. It has then three important properties that can be
summarized in the next theorem (see
\cite{Nor}).

\begin{Theorem} \label{Northcottideal}
Let ${\bf v}$, ${\bf u}$ and
$\varphi$ as above. Then either $ N $ is the whole ring or it is a
proper ideal with the following properties
\begin{itemize}
\item[{\rm (a)}] $N$ is  a perfect ideal (i.e. its grade is equal to
the projective dimension of $ R/N$);
\item[{\rm (b)}] $ ({\bf v}) \colon N  = ({\bf u}) $;
\item[{\rm (c)}] $ ({\bf v})\colon ({\bf u}) = N$.
\end{itemize}
\end{Theorem}

The case where $({\bf u})$  is the maximal ideal of a regular local ring is
particularly interesting: $\det \varphi$ will then generate the socle
of the Artin algebra $R/({\bf v})$.

\bigskip

\noindent {\bf Proof of Theorem~\ref{primelink1}.}
According to Corollary~\ref{C1}, we must show that $\delta(I)\neq 0$,
as this leads to $\lambda(I^2/JI)=0$. As observed earlier,
$\delta(I)=0$
 means that $I$ is also generated by a regular
sequence. This puts us in the context of Northcott ideals. Since
however $\lambda(I/J)=1$, it follows that the ideal $N = (J, \det
\varphi)$ is perfect and necessarily the maximal ideal of $R$. In the
case that $R$ is not regular, this is a contradiction. In the other
case, when at least two of the generators of $J$ lie in the square of
the maximal ideal ${\frak m}$, we cannot have the equality ${\frak m}=(J,
\det\varphi )$ either. \QED

\section{Cohen--Macaulay Rees Algebras}

In this section  we are going to use the results obtained in the previous
section---in particular Theorem~\ref{primelink1}---in order to study the
properties of the Rees algebras of the links of  Cohen--Macaulay prime
ideals. (We use the terminology and some of the methods described in
\cite{HSV1} and some of its references.)

\begin{Theorem}\label{main}
Let $ R $ be a $ d$-dimensional Gorenstein local ring,
$ {\frak p} $ be a Cohen--Macaulay prime ideal of codimension $g$,
$J= (a_1,...,a_{g})$ a complete intersection, and let $I=J\colon{\frak p}$.
If one of the following holds:
\begin{itemize}
\item[{\rm (a)}]
 $R_{{\frak p}}$ is not a regular local ring, or
\item[{\rm (b)}]
$R_{{\frak p}}$ is a regular local ring and  two of the generators of $J$
belong  to ${\frak p}^{(2)}$,
\end{itemize}
then the associated graded ring $ {\rm gr}_I(R) $ is Cohen--Macaulay.
\end{Theorem}

\demo First observe that $I$ is a Cohen--Macaulay ideal
(see \cite{pszpiro}).
Let $R[Jt]$ and $R[It]$ be the Rees algebras of the ideals $J$ and
$I$. Since $I^2=JI$, by Theorem~\ref{primelink1}, $R[It]$ is a
finitely generated $R[Jt]$--modules and  $I\cdot R[Jt] = I\cdot
R[It]$.

On the other hand, the exact sequence
$$
\diagram
0 \rto & I\cdot R[Jt] \rto & R[Jt] \rto & R[Jt]\otimes_R R/I \rto & 0
\enddiagram
$$
implies that $I\cdot R[Jt]$ is a maximal Cohen--Macaulay module over $R[Jt]$,
since $R[Jt]$ is a Cohen--Macaulay ring of dimension $ d+1 $ and
the fibre $R[Jt]\otimes_R R/I$ is a polynomial ring in $ g $
variables over the Cohen--Macaulay ring
$R/I$; it is thus a Cohen--Macaulay ring of dimension $\dim(R/I) +
g=d-g+g=d$.   This means that, when viewed as an
$R[It]$--module, $I\cdot R[It]$ is a maximal Cohen--Macaulay module as
well.

The assertion about the associated graded ring of $I$,
$$
{\rm gr}_I(R)=\sum_{n\ge 0} I^{n}/I^{n +1} = R[It]\otimes_R R/I
$$
will follow from depth counting in the
 exact sequences:
\begin{equation}\label{hunekeeq1}
\diagram
0 \rto & I\cdot R[It] \rto & R[It] \rto & {\rm gr}_I(R) \rto & 0,
\enddiagram
\end{equation}
\begin{equation}\label{hunekeeq2}
\diagram
0 \rto & It\cdot R[It] \rto & R[It] \rto & R \rto & 0,
\enddiagram
\end{equation}
and the canonical isomorphism $It\cdot R[It]\simeq I\cdot R[It]$. \QED

\begin{Corollary} If in {\rm Theorem~\ref{main}} $ d \ge 2 $
$($which is automatic in the second case$)$, then the Rees algebra of
$I$ is Cohen--Macaulay.
\end{Corollary}

\demo Since the reduction number of $I$ is 1, the assertion follows
from the Theorem and \cite{IT}. \QED

\begin{Remark}{\rm
Although Theorem~\ref{main} covers most of the cases of links of prime
ideals, the restriction in (b) leaves interesting links out. For
example, suppose ${\frak p}$ is a Gorenstein prime ideal of
codimension $ g $ and let $ J $ be a
regular sequence in ${\frak p}$ with $ g $ elements. The link
$I = J\colon {\frak p}$ is generated by $J$ with an extra element, $I = (J,
z)$. If $I_{{\frak p}}= {\frak p} R_{{\frak p}}$, as it may
as well occur, then obviously
we are not going to have $I^2=JI$. In fact, in this case, $I$ is going
to be generated by analytically independent elements and therefore it
will not admit a proper reduction. Nevertheless, one can still
show that the Rees algebra of $I$ is Cohen--Macaulay (see
\cite{HSV1}).

The authors are also aware of other ideals, in non Gorenstein rings,
with some of these properties. But they arise in much more
constrainted settings.
}\end{Remark}

\end{document}